\documentclass[runningheads]{llncs}
\usepackage[utf8]{inputenc}
\usepackage[T1]{fontenc}

\usepackage{graphicx}
\usepackage{mathtools}
\mathtoolsset{showonlyrefs}
\usepackage{color}

\usepackage[hyphens]{url}   
\usepackage[unicode]{hyperref} 

\usepackage{amsthm}

\usepackage[vvarbb]{newtxmath}

\usepackage{algorithm}
\usepackage{algorithmic}
\usepackage{enumitem}
\usepackage{subcaption}
\usepackage{etoolbox}

\usepackage[vvarbb]{newtxmath}
\usepackage{newtxtext}

\floatstyle{ruled}
\newfloat{oracle}{htbp}{ora}
\floatname{oracle}{Oracle}

\theoremstyle{definition}
\theoremstyle{remark}

\DeclarePairedDelimiterX{\normSimple}[1]{\lVert}{\rVert}
{\ifblank{#1}{\mathord{\cdot}}{#1}}

\newcommand{\norm}[2][]{\normSimple{#2}\ifblank{#1}{}{\sb{#1}}}

\DeclarePairedDelimiter{\card}{\lvert}{\rvert}
\DeclarePairedDelimiterXPP{\supp}[1]{\operatorname{supp}}(){}{#1}

\newcommand{\R}{\mathbb R}
\newcommand{\Z}{\mathbb Z}

\newcommand{\innp}[2]{\left\langle #1, #2 \right\rangle}

\newcommand*{\set}[2]{\left\{{#1}\,\middle|\,{#2}\right\}}

\DeclareMathOperator*{\conv}{conv}

\DeclareMathOperator*{\argmax}{argmax}
\DeclareMathOperator*{\argmin}{argmin}

\DeclarePairedDelimiter\ceil{\lceil}{\rceil}
\DeclarePairedDelimiter\floor{\lfloor}{\rfloor}

\begin{document}
\title{Restarting Algorithms: \\ Sometimes there is Free Lunch}
\titlerunning{Restarting Algorithms}
%
\author{Sebastian Pokutta\inst{1}\inst{2}}
\authorrunning{S. Pokutta}
%
\institute{Technische Universität Berlin \and
Zuse Institute Berlin \\
\email{pokutta@zib.de}
}
\maketitle              
\begin{abstract}
 In this overview article we will consider the deliberate restarting of
 algorithms, a meta technique, 
 in order to improve the algorithm's performance, e.g., convergence
 rates or approximation guarantees. One of the major advantages is that
 restarts are \emph{relatively} black box, not requiring any
 (significant) changes to the base algorithm that is
 restarted or the underlying argument, while leading to \emph{potentially significant}
 improvements, e.g., from sublinear to linear rates of
 convergence. Restarts are widely used in different fields and have
 become a powerful tool to leverage additional information that has
 not been directly incorporated in the base algorithm or argument. We
 will review restarts in various settings from continuous
 optimization, discrete optimization, and submodular function
 maximization where they have delivered impressive results. 

\end{abstract}

\section{Introduction}

Restarts are a powerful meta technique to improve the behavior of
algorithms. The basic idea is to deliberately restart some base algorithm,
often with changed input parameters, to speed-up convergence, improve
approximation guarantees, reduce number of calls to expensive
subroutines and many more, often leading to provably
better guarantees as well as significantly improved real-world
computational performance. In actuality this comes down to running a
given algorithm with a given set of inputs for some number of
iterations, then changing the input parameters usually as a function of
the output, and finally restarting the algorithm with new input
parameters; rinse and repeat. 

One appealing aspect of restarts is that they are relatively
black-box, requiring only little to no knowledge of the to-be-restarted base
algorithm except for the guarantee of the base algorithm that is then
amplified by means of restarts. The reason why restarts often work, i.e.,
improve the behavior of the base algorithm is that some
structural property of the problem under consideration is not
explicitly catered for in the base algorithm, e.g., the base algorithm
might work for general convex functions, however the function under
consideration might be strongly convex or sharp. Restarts cater to
this additional problem structure and are in particular useful
when we want to incorporate data-dependent parameters. In fact, for
several cases of interest the only known way to incorporate that
additional structure is via restarts often pointing out a missing piece
in our understanding.

On the downside, restarts often explicitly depend on parameters
arising from the additional structure under consideration and obtained guarantees are off
by some constant factor or even log factor. The former can be often
remedied with adaptive or scheduled restarts (see e.g.,
\cite{Roul17,hinder2020generic}) albeit with some minor cost. This way
we can obtain fully adaptive algorithms that adapt to additional
structure without knowing the accompanying parameters explicitly. The
latter shortcoming is inherent to restart scheme as due to their black
box nature additional structural information might not be incorporated
perfectly.

Restarts have been widely used in many areas and fields and we will
review some of these applications below to provide context. We would
like to stress that references will be incomplete and biased; please
refer also to the references cited therein.

\paragraph{SAT solving and Constraint Programming.} Restarts are
ubiquitous in SAT Solving and Constraint Programming to, e.g., explore
different parts of the search space. Also after new clauses have been
learned, these clauses are often added back to the formulation and
then the solver is restarted. This can lead to dramatic overall
performance improvements for practical solving; see e.g.,
\cite{huang2007effect,biere2008adaptive} and references contained
therein.

\paragraph{Global Optimization.} Another important area where restarts
are used is global optimization. Often applied to non-convex problems,
the hope is that with randomized restarts different local optima can
be explored, ideally one of those being a global one; see e.g.,
\cite{hu2009random} and their references.

\paragraph{Integer Programming.} Modern integer programming solvers
use restarts in many different ways, several of which have been
inspired by SAT solving and Constraint Programming. In fact, Integer
Programming solvers can be quite competitive for pseudo-Boolean
problems \cite{berthold2008solving}. A relatively recent approach
\cite{anderson2019clairvoyant} is \emph{clairvoyant restarts} based on
online tree-size estimation that can significantly improve solving
behavior. \medskip

Most of the restart techniques mentioned above, while very important,
come without strong guarantees. In this article, we are more
interested in cases, where provably strong guarantees can be obtained
that also translate into real-world computational advantages. In the
following, we will restrict the discussion to three examples from
convex optimization, discrete optimization, and submodular function
maximization. However, before we consider those, we would like to
mention a two related areas where restarts have had a great impact not
just from a computational point of view but also to establish new
theoretical guarantees, but that are unfortunately beyond the scope of
this overview.
 
\paragraph{Variance Reduction via Restarts.} Usually when we consider
stochastic convex optimization problems where the function is given as
a general expectation and we would like to use first-order methods for
solving the stochastic problem, we cannot expect a convergence rate
better than $O(1/\sqrt{t})$ under usual assumptions, where $t$ is the
number of stochastic gradient evaluations. However, it turns out that
if we consider so-called finite sum problems, a problem class quite
common in machine learning, where the expectation is actually a finite
sum and some mild additional assumptions are satisfied, then we can
obtain a linear rate of convergence by means of variance
reduction. This is an exponential improvement in convergence
rate. Variance reduction techniques replace the stochastic gradient
which is an unbiased estimator of the true gradient with a different,
lower variance, unbiased estimator that is formed with the help of a
reference point obtained from an earlier iterate. This reference point
is then periodically reset via a restart scheme. Important algorithms
here are for example Stochastic Variance Reduced Gradient Descent
(SVRG) \cite{johnson2013accelerating} and its numerous variants, such
as e.g., the Stochastic Variance Reduced Frank-Wolfe algorithm (SVRFW)
\cite{hazan2016variance}.

\paragraph{Acceleration in Convex Optimization.} Restarts have been
heavily used in convex optimization both for improving convergence
behavior via restarts in real-world computations (see e.g.,
\cite{o2015adaptive}) but also as part of formal arguments to
establish accelerated convergence rates and design provably faster
algorithms. As the literature is particularly deep, we will sample only a few
of those works in the context of first-order methods here that we are
particularly familiar with; we provide further references in the
sections to come. For example restarts have been used in
\cite{allen2014linear} to provide an alternative explanation of
Nesterov's acceleration as arising from the coupling mirror descent
and gradient descent. In \cite{Roul17} it has been shown how restarts
can be leveraged to obtain improved rates as the sharpness of the
function (roughly speaking how fast the function curves around its
minima) increases and these restart schemes have been also
successfully carried over to the conditional gradients case in
\cite{KDP2018}. Restarts have been also used to establish
dimension-independent local acceleration for conditional gradients
\cite{diakonikolas2019locally} by means of coupling the Away-step
Frank-Wolfe algorithm with an accelerated method. As we will see later
in the context of submodular maximization, restarts can be also used
to reduce the number of calls to expensive oracles. This have been
extensively used for lazification of otherwise expensive algorithms in
\cite{BPZ2017,BPZ2017jour} leading to several orders of speed-up in
actual computations while maintaining worst-case guarantees identical
to those of the original algorithms and in \cite{LPZZ2017} a so-called
optimal method based on lazification has been derived. Very recently,
in \cite{hinder2020generic} a new adaptive restart scheme has been
presented that does not require any knowledge of otherwise
inaccessible parameters and its efficacy for saddle point problems has
been demonstrated.

\subsection*{Outline} In Section~\ref{sec:smoothopt} we consider
restart examples from convex optimization and in
Section~\ref{sec:discopt} we consider examples from discrete
optimization. Finally we consider submodular function maximization in
Section~\ref{sec:sfm}. We keep technicalities to a bare minimum,
sometimes simplifying arguments for the sake of exposition. We provide
references though with the full argument, for the interested reader.

\section{Smooth Convex Optimization}
\label{sec:smoothopt}

Our first examples come from smooth convex optimization. As often, the
examples here are (arguably) the cleanest ones. We briefly recall
some basic notions:

\begin{definition}[Convexity]
  Let \(f: \R^n \rightarrow \R\) be a differentiable function. Then
  \(f\) is \emph{convex}, if for all \(x,y\) it holds:
  \begin{equation}
    \label{eq:convex}
      f(y) - f(x) \geq \innp{\nabla f(x)}{y-x}.
    \end{equation}
    In particular, all local mimima of \(f\) are global minima of
    \(f\). 
\end{definition}

\begin{definition}[Strong Convexity]
  \label{def:SC}
  Let \(f: \R^n \rightarrow \R\) be a differentiable convex function. Then
  \(f\) is \emph{\(\mu\)-strongly convex} (with \(\mu > 0\)), if for all \(x,y\) it holds:
  \begin{equation}
    \label{eq:Sconvex}
      f(y) - f(x) \geq \innp{\nabla f(x)}{y-x} + \frac{\mu}{2} \norm{y-x}^2.
    \end{equation}
\end{definition}

\begin{definition}[Smoothness]
  Let \(f: \R^n \rightarrow \R\) be a differentiable function. Then
  \(f\) is \emph{\(L\)-smooth} (with \(L >0\)), if for all \(x,y\) it holds:
  \begin{equation}
    \label{eq:smooth}
      f(y) - f(x) \leq \innp{\nabla f(x)}{y-x} + \frac{L}{2} \norm{y-x}^2.
    \end{equation}
\end{definition}

In the following let \(x^* \in X^* = \argmin f(x)\) denote an optimal
solution from the set of optimal solutions \(X^*\). Choosing
\(x = x^*\) and applying the definition of strong convexity (Definition~\ref{def:SC}) we
immediately obtain:

\begin{equation}
  \label{eq:primalbound}
  f(y) - f(x^*) \geq \innp{\nabla f(x^*)}{y-x^*} + \frac{\mu}{2}
  \norm{y-x^*}^2 \geq \frac{\mu}{2}
  \norm{y-x^*}^2,
\end{equation}
where the last inequality follows from \(\innp{\nabla
  f(x^*)}{y-x^*}\geq 0\) by first-order optimality of \(x^*\) for \(\min f(x)\), i.e., the primal
gap upper bounds the distance to the optimal solution. This also
implies that the optimal solution \(x^*\) is unique. 

\subsubsection*{Smooth Convex to Smooth Strongly Convex: the basic case.}

Let \(f: \R^n \rightarrow \R\) be an $L$-smooth convex function. Then using
\emph{gradient descent}, updating iterates \(x_t\) according to \(x_{t+1} \leftarrow x_t - \frac{1}{L}
\nabla f(x_t)\), yields the following standard guarantee, see e.g., \cite{nemirovski2001lectures,nesterov2018introductory,hazan2019lecture,lanfirst}

\begin{proposition}[Convergence of gradient descent: smooth convex
  case]
  \label{prop:gd-conv}
  Let \(f: \R^n \rightarrow \R\) be a smooth convex function and \(x_0
  \in \R^n\) and \(x^* \in X^*\). Then gradient descent generates a sequence of iterates satisfying
  \begin{equation}
    \label{eq:gd-conv}
    f(x_t) - f(x^*) \leq \frac{L\norm{x_0 - x^*}^2}{t}.
  \end{equation}
\end{proposition}

Now suppose we additionally know that the function \(f\) is
\(\mu\)-strongly convex. Usually, we would expect a linear rate of
convergence in this case, i.e., to reach an additive error of
\(\varepsilon\), we would need at most
\(T \leq \frac{L}{\mu} \log \frac{f(x_0) - f(x^*)}{\varepsilon}\)
iterations. However, rather than reproving the convergence rate (which
is quite straightforward in this case) we want to reuse the guarantee
in Proposition~\ref{prop:gd-conv} as a black box and the \(\mu\)-strong convexity of
\(f\). We will use the simple restart scheme given in
Algorithm~\ref{alg:simple-restart}: in restart phase \(\ell\) we run a given base
algorithm \(\mathcal A\) for a fixed number of iterations \(T_\ell\)
on the iterate \(x^{\ell-1}\) output in the previous iteration:

\begin{algorithm}[htbp]
\caption{Simple restart scheme}
\label{alg:simple-restart}
\begin{algorithmic}[1]
\REQUIRE Initial point $x_0 \in \R^n$, base algorithm \(\mathcal A\),
iteration counts \((T_\ell)\).
\ENSURE Iterates $x^1, \dotsc, x^K \in \R^n$.
\FOR{$\ell=1$ \TO $K$} 
  \STATE $x^\ell \leftarrow \mathcal A(f,x^{\ell-1},T_\ell)$
  \COMMENT{run base algorithm for \(T_\ell\) iterations} 
\ENDFOR
\end{algorithmic}
\end{algorithm}

A priori, it is unclear whether the restart scheme in
Algorithm~\ref{alg:simple-restart} is doing anything useful, in fact
even convergence might not be immediate as we in principle could undo
work that we did in a preceding restart phase. Also note that when
restarting vanilla gradient descent with a fixed step size of
\(\frac{1}{L}\) as we do here the final restarted algorithm is
identical to vanilla gradient descent, i.e., the restarts do not
change the base algorithm. This might seem nonsensical and
we will get back to this soon; the reader can safely ignore this for now. 

In order to analyze our restart scheme we
first chain together Inequalities~\eqref{eq:gd-conv} and
\eqref{eq:primalbound} and obtain:

\begin{equation}
  \label{eq:gd-restart-inequality}
  f(x_t) - f(x^*) \leq \frac{L\norm{x_0 - x^*}^2}{t} \leq
    \frac{2L}{\mu} \ \frac{f(x_0) - f(x^*)}{t}. 
\end{equation}

This chaining together of two error bounds is at the core of most
restart arguments and we will see several variants of this. Next we
estimate how long we need to run the base method, using
Inequality~\eqref{eq:gd-restart-inequality} to halve the primal gap
from some given starting point $x_0$ (this will be the point from
which we are going to restart the base method), i.e., we want to find
\(t\) such that

\begin{equation}
  \label{eq:gd-time-to-halve}
  f(x_t) - f(x^*)\leq
    \frac{2L}{\mu} \ \frac{f(x_0) - f(x^*)}{t} \leq \frac{f(x_0) - f(x^*)}{2},
  \end{equation}
which implies that it suffices to run gradient descent for \(T_\ell \doteq
  \left \lceil
  \frac{4L}{\mu}\right \rceil\) steps for all $\ell = 1, \dots, K$ to
halve a given primal bound as there is no dependency on the state of
the algorithm in
this case. Now, in
order to reach \(f(x_T) - f(x^*) \leq \varepsilon\), we have to halve
\(f(x_0) - f(x^*)\) at most \(K \doteq \left \lceil \log \frac{f(x_0) - f(x^*)}{\varepsilon}
\right \rceil\) times and each of the halving can be accomplished in
at most \(\left \lceil \frac{4L}{\mu}\right \rceil\) gradient descent
steps. All in all we obtain that after at most
\begin{equation}
  \label{eq:gd-restart-rate}
  T \geq \sum_{\ell = 1}^{K} T_\ell = K \cdot T_1 = \left \lceil \frac{4L}{\mu}\right \rceil \left \lceil \log \frac{f(x_0) - f(x^*)}{\varepsilon}
\right \rceil
\end{equation}
gradient descent steps we have obtained a solution
\(f(x_T) - f(x^*) = f(x^K) - f(x^*)\leq \varepsilon\). With this we
have obtained the desired convergence rate. Note that the iterate
bound in Inequality~\eqref{eq:gd-restart-rate} is optimal for vanilla
gradient descent up to a constant factor of $4$; see e.g.,
\cite{nemirovski2001lectures,lanfirst,hazan2019lecture}.

In the particular case from above it is also important to observe that
our base algorithm gradient descent is essentially memoryless. In
fact, the restarts do not \lq{}reset\rq{} anything in this particular
case and so we have also indirectly proven that gradient descent
\emph{without} restarts will converge with the rate from
Inequality~\eqref{eq:gd-restart-rate}. This is particular to this
example though and will be different in our next one. Also, note that
a direct estimation would have yielded the same rate up to the
factor \(4\) discussed above.

\subsubsection{Smooth Convex to Smooth Strongly Convex: the
  accelerated case.}
While the rate from Inequality~\eqref{eq:gd-restart-rate} is
essentially optimal for vanilla gradient descent it is known that
(vanilla) gradient descent itself is not optimal for smooth and
strongly convex functions and also Proposition~\ref{prop:gd-conv} is not
optimal for smooth and (non-strongly) convex functions. In fact
Nesterov showed in \cite{nesterov1983method} that for smooth and
(non-strongly) convex functions a quadratic improvement can be obtained;
a phenomenon commonly referred to as \emph{acceleration}:

\begin{proposition}[Convergence of accelerated gradient descent]
  \label{prop:gd-conv-acc}
  Let \(f: \R^n \rightarrow \R\) be an $L$-smooth convex function and \(x_0
  \in \R^n\) and \(x^* \in X^*\). Then accelerated gradient descent generates a sequence of iterates satisfying
  \begin{equation}
    \label{eq:gd-conv-acc}
    f(x_t) - f(x^*) \leq \frac{cL\norm{x_0 - x^*}^2}{t^2},
  \end{equation}
  for some constant \(c > 0\). 
\end{proposition}

Again, we could try to directly prove a better rate via acceleration
for the smooth and strongly case (which is non-trivial this time) or,
as before, invoke our restart scheme in
Algorithm~\ref{alg:simple-restart} in a black-box fashion, which is
what we will do here. As before we will use an analog of
Inequality~\eqref{eq:gd-restart-inequality} to estimate how long it
takes to halve the primal gap, i.e., we want to find
\(t\) such that

\begin{equation}
  \label{eq:agd-time-to-halve}
  f(x_t) - f(x^*)\leq
    \frac{2cL}{\mu} \ \frac{f(x_0) - f(x^*)}{t^2} \leq \frac{f(x_0) - f(x^*)}{2},
  \end{equation}
  which implies that it suffices to run accelerated gradient descent
  for
  \(T_\ell \doteq \left \lceil \sqrt{\frac{4cL}{\mu}}\right \rceil\)
  steps for all $\ell = 1, \dots, K$ to halve a given primal gap. With
  the same reasoning as above we need to halve the primal gap at most
  \(K \doteq \left \lceil \log \frac{f(x_0) - f(x^*)}{\varepsilon}
  \right \rceil\) times to reach an additive error of $\varepsilon$. Putting
  everything together we obtain that after at most

\begin{equation}
  \label{eq:agd-restart-rate}
  T \geq \sum_{\ell = 1}^{K} T_\ell = K \cdot T_1 = \left \lceil \sqrt{\frac{4cL}{\mu}}\right \rceil \left \lceil \log \frac{f(x_0) - f(x^*)}{\varepsilon}
\right \rceil
\end{equation}
accelerated gradient descent steps we have obtained a solution
\(f(x_T) - f(x^*) = f(x^K) - f(x^*)\leq \varepsilon\). Note that the
iterate bound in Inequality~\eqref{eq:agd-restart-rate} is optimal for
strongly convex and smooth functions (up to a constant factor). In
contrast to the unaccelerated case, this time the restart actually \lq{}resets\rq{} the base
algorithm as accelerated gradient descent uses a specific step size
strategy that is then reset.

\begin{remark}
  Sometimes it is also possible to go backwards. Here we recover the
  optimal base algorithm for the smooth and (non-strongly) convex case
  from the strongly convex one. The argument is due to
  \cite{nesterov2012make} (we follow the variant in
  \cite{scieur2016regularized}). Suppose we know an optimal algorithm
  for the strongly convex and smooth case that ensures
  $f(x_T) - f(x^*) \leq \varepsilon$ after
  $O\left(\sqrt{\frac{L}{\mu}} \log \frac{f(x_0) -
      f(x^*)}{\varepsilon}\right)$ iterations. Now consider a smooth and convex function $f$ and an initial iterate $x_0$
  together with some upper bound $D$ on the distance to some optimal
  solution, i.e., $\norm{x_0 - x^*} \leq D$. Given an accuracy
  $\varepsilon > 0$, we consider the auxiliary function
$$f_\varepsilon(x) \doteq f(x) + \frac{\varepsilon}{2D^2} \norm{x-x_0}^2,$$
which is $\left(L + \frac{\varepsilon}{2D^2}\right)$-smooth and
$\frac{\varepsilon}{2D^2}$-strongly convex. It can be easily seen that
$$f(x) - f(x^*) \leq f_\varepsilon(x) - f_\varepsilon(x^*) + \frac{\varepsilon}{2},$$
so that finding an $\varepsilon/2$-optimal solution to $\min f_\varepsilon$ provides an $\varepsilon$-optimal solution to $\min f$. We can now run the purported optimal method on the smooth and strongly convex function $f_\varepsilon$ to compute an $\varepsilon/2$-optimal solution to $\min f_\varepsilon$, which we obtain after:
$$
O\left(\sqrt{\frac{L + \frac{\varepsilon}{2D^2}}{\frac{\varepsilon}{2D^2}}} \log 2 \frac{f_\varepsilon(x_0) - f_\varepsilon(x^*)}{\varepsilon} \right) \leq O\left(\sqrt{\frac{2LD^2 + \varepsilon}{\varepsilon}} \log \frac{(L+\varepsilon)D^2}{\varepsilon} \right),
$$
iterations, where we used $f_\varepsilon(x_0) - f_\varepsilon(x^*) \leq \frac{(L+\varepsilon)D^2}{2}$. Finally note, ignoring the log factor, $\sqrt{\frac{2LD^2 + \varepsilon}{\varepsilon}} \leq T \Leftrightarrow \frac{2LD^2 + \varepsilon}{T^2} \leq \varepsilon$, which is the bound from Proposition~\ref{prop:gd-conv-acc}. 
\end{remark}

The approach used in this section to obtain better rates of convergence under
stronger assumptions by means of the simple restart scheme in
Algorithm~\ref{alg:simple-restart} works in much broader settings in
convex optimization (including the constrained case). For example it
can be used to improve the $O(1/\sqrt{t})$-rate for general non-smooth
convex functions via sub-gradient descent into the $O(1/t)$-rate for
the non-smooth strongly convex case. Here the base rate is
$f(x_t) - f(x^*) \leq \frac{G \norm{x_0-x^*}}{\sqrt{t}}$, where $G$ is a
bound on the norm of the subgradients. We obtain the restart
inequality chain (analog to
Inequality~\eqref{eq:gd-restart-inequality}):

\begin{equation}
  \label{eq:subgd-restart-inequality}
  f(x_t) - f(x^*) \leq \frac{G \norm{x_0-x^*}}{\sqrt{t}} \leq
    \frac{G}{\sqrt{t}} \ \sqrt{\frac{f(x_0) - f(x^*)}{\mu}},
\end{equation}
and halving the primal gap takes at most $\frac{4G^2}{\mu (f(x_0) - f(x^*))}$ iterations. Following the argumentation from above, we then arrive that the total number of required subgradient descent iterations using Algorithm~\ref{alg:simple-restart} to ensure $f(x_t) - f(x^*) \leq \varepsilon$ is at most $t \geq \frac{8G^2}{\varepsilon \mu}$ for the non-smooth but $\mu$-strongly convex case, which is optimal up to constant factors.  

\paragraph{Related approaches.}

In a similar way we can incorporate additional information obtained
e.g., from so-called \emph{Hölder(ian) Error Bounds} or
\emph{sharpness} (see, e.g.,
\cite{bolte2007lojasiewicz,bolte2017error} and references contained
therein for an overview). The careful reader might have observed that
the restart scheme in Algorithm~\ref{alg:simple-restart} requires
knowledge of the parameter $\mu$. While this could be acceptable in
the strongly convex case, for more complex schemes to leverage, e.g.,
sharpness, this is unacceptable as the required parameters are hard to
estimate and generally inaccessible. This however, can be remedied in
the case of sharpness, at the cost of an extra $O(\log^2)$-factor in the
rates, via \emph{scheduled restarts} as done in \cite{Roul17} that do
not require sharpness parameters as input or when an error bound (of
similar convergence rate) is available as in the case of conditional
gradients \cite{KDP2018}; see also \cite{hinder2020generic} for a very
recent adaptive restart scheme using error bounds estimators. 

\section{Discrete Optimization}
\label{sec:discopt}

In this section we consider a prominent example from integer
programming: optimization via augmentation, i.e., optimizing by
iteratively improving the current solution. 

We consider the problem:

\begin{equation}
  \label{eq:ip}
  \max \set{cx}{x \in P \cap \Z^n},
\end{equation}
where \(P \subseteq \R^n\) is a polytope and \(c \in \Z^n\).

\begin{algorithm}[htbp]
 \caption{\label{alg:augmentation}Augmentation} 
\begin{algorithmic}[1]
\REQUIRE Feasible solution \(x^0\) and objective \(c \in \Z_+^n\)\\
\ENSURE Optimal solution of \(\max \set{cx}{x \in P \cap \{0,1\}^n}\)
    \STATE \(\tilde{x} \leftarrow x^0\)
    \REPEAT
    \STATE \label{eq:augStep}\textbf{compute} \(x \in P\) integral with \(c(x -
    \tilde{x}) > 0\) and set \(\tilde{x} \leftarrow x\)
    \COMMENT{improve solution} 
    \UNTIL{no improving solution exists}
  \RETURN \(\tilde{x}\) \COMMENT{return optimal solution}
  \end{algorithmic}
\end{algorithm}

To simplify the exposition we
assume that \(P \subseteq [0,1]^n\) and \(c \geq 0\) (the latter is
without loss of generality by flipping coordinates), however the arguments here
generalize to the general integer programming case. Suppose further
that we can compute \emph{improving solutions}, i.e., given \(c\) and
a solution \(x_0\), we can compute a new solution \(x\), so
that \(c(x - x_0) > 0\) if \(x_0\) was not already optimal; such a step (Line~\ref{eq:augStep} in Algorithm~\ref{alg:augmentation}) is called an \emph{augmentation step}. Then a
trivial and inefficient strategy is Algorithm~\ref{alg:augmentation},
where we continue improving the solution until we have reached the
optimum.  It is not too hard to see that Algorithm~\ref{alg:augmentation} can
take up to \(2^n\) steps, essentially enumerating all feasible
solutions to reach the optimal solution; simply consider the cube \(P
= [0,1]^n\) and an objective \(c\) with powers of \(2\) as entries.

\subsubsection{Bit Scaling.}
\label{sec:bitscaling}

We will now show that we can do significantly better by restarting
Algorithm~\ref{alg:augmentation}, so that we obtain a number of augmentation steps
of \(O(n \log \norm{c}_\infty)\), where
\(\norm{c}_\infty \doteq \max_{i \in [n]} c_i\). This is an
exponential improvement over base algorithm and the restart scheme,
called \emph{bit scaling}, is due to \cite{schulz19950} (see also
\cite{edmonds1972theoretical,graham1995handbook}). It crucially relies
on the following insight: Suppose we decompose our objective \(c = 2c_1 +
c_0\) with \(c_0 \in \{0,1\}^n\) (note this decomposition is unique) and we have already obtained some
solution \(x_0 \in P \cap \{0,1\}^n\) that is optimal for \(\max \set{c_1x}{x \in P \cap
  \Z^n}\), then we have for all \(x \in P \cap \{0,1\}^n\):
\begin{equation}
  \label{eq:bit-restart}
  c(x-x_0) = 2\underbrace{c_1(x-x_0)}_{\leq 0} + c_0 (x-x_0) \leq n,
\end{equation}
by the optimality of \(x_0\) for \(c_1\) and
\(c_0, x, x_0 \in \{0,1\}^n\). Hence starting from \(x_0\), for
objective \(c\), there are at most \(n\) augmentation steps to be
performed with Algorithm~\ref{alg:augmentation} to reach an optimal
solution for \(c\). Equipped with Inequality~\eqref{eq:bit-restart} the
following strategy emerges: slice by the objective \(c\)
according to its bit representation and then successively optimize
with respect to the starting point from a previous slice. We first
present the formal bit scaling restart scheme in
Algorithm~\ref{alg:bitScaling}, where \(\mathcal A\) denotes
Algorithm~\ref{alg:augmentation}.

\begin{algorithm}[htbp]
 \caption{\label{alg:bitScaling}Bit Scaling} 

\begin{algorithmic}[1]
\REQUIRE Feasible solution \(x^0\)\\
\ENSURE Optimal solution to \(\max \set{cx}{x \in P \cap \Z^n}\)
    \STATE \(C \leftarrow \norm{c}_\infty + 1\), \(\mu \leftarrow 2^{\ceil{\log C}}\), \(\tilde{x}
    \leftarrow x^0\), \(c^\mu \leftarrow \floor{c / \mu}\) \COMMENT{initialization}
    \REPEAT
    \STATE \textbf{Call} \(\tilde x \leftarrow \mathcal A(\tilde x, c^\mu)\) \label{aeq:callStdAug}
    \STATE \(\mu \leftarrow \mu/2\), \(c^\mu \leftarrow \floor{c / \mu}\)
    \UNTIL{\(\mu < 1\)} \label{aeq:stop}
  \RETURN \(\tilde{x}\) \COMMENT{return optimal solution}
  \end{algorithmic}
\end{algorithm}

Next, we will show that restart scheme from
Algorithm~\ref{alg:bitScaling} requires at most
$O(n \log \norm{c}_\infty)$ augmentation steps (Line~\ref{eq:augStep}
in Algorithm~\ref{alg:augmentation}) to solve
Problem~\eqref{eq:ip}. First observe, that by construction and the
stopping criterion in Line~\ref{aeq:stop} of
Algorithm~\ref{alg:bitScaling} it is clear that we call $\mathcal A$
in Line~\ref{aeq:callStdAug} at most $\lceil \log C \rceil$
times. Next, we bound the number of augmentation steps in
Line~\ref{aeq:callStdAug}  executed
within algorithm $\mathcal A$. To this end, let $\tilde x$
and $\mu$ denote the input to $\mathcal A$. In the first iteration
$c^\mu \in \{0,1\}^n$, so that $\mathcal A$ can perform at most $n$
augmentation steps. For later iterations observe that $\tilde x$ was
optimal for $c^{2\mu} = \lfloor c / (2\mu) \rfloor$. Moreover, we have
$c^{\mu} = \lfloor c / \mu \rfloor = 2 c^{2\mu} + c_0$, where
$c_0 \in \{0,1\}^n$ as before. Via
Inequality~\eqref{eq:bit-restart} we obtain for all feasible solutions
\(x \in P \cap \Z^n\):
\begin{equation}
  c^{\mu}(x-\tilde x) = 2c^{2\mu}(x-\tilde x) + c_0 (x-\tilde x) \leq n,	
\end{equation}
which holds in particular for the optimal solution $x^*$ to Problem~\eqref{eq:ip}. As each augmentation step reduces the primal gap $c^{\mu}(x-\tilde x)$ by at least $1$, we can perform at most $n$ augmentation steps. This completes the argument.

\subsubsection{Geometric Scaling.} The restart scheme in
Algorithm~\ref{alg:bitScaling} essentially restarted via bit-scaling the
objective function, hence the name. We will now present a more versatile restart
scheme that is due to \cite{schulz2002complexity} (see also \cite{BPPP2015}
for a comparison and worst-case examples),
which essentially works by restarting a regularization of our
objective \(c\). For
comparability we also consider Problem~\eqref{eq:ip} here, however the
approach is much more general, e.g., allowing for general integer
programming problems and with modifications even convex
programming problems
over integers. 

Again, we will modify the considered objective function $c$ in each restart. Given the original linear objective $c$, we will consider:

\begin{equation}
c^{\mu}(x,\tilde x) = c(x - \tilde x) - \mu \norm{x - \tilde x}_1. 
\end{equation}

Note that $c^{\mu}(x, \tilde x)$ is a linear function in $x \in \{0,1\}^n$ for a
given $\tilde x \in \{0,1\}^n$.
In particular we can call
Algorithm~\ref{alg:augmentation} with objective
$c^{\mu}(\cdot, \cdot)$ and starting point $\tilde x$. The restart
scheme works as follows: For a given $\mu$ we call
Algorithm~\ref{alg:augmentation} with objective
$c^{\mu}(\cdot, \cdot)$ and starting point $\tilde x$. Then we
halve $\mu$ and repeat.

As in the bit-scaling case, the key is to estimate the number of augmentation steps
performed in such a call. To this end let $x_0$ be returned by
Algorithm~\ref{alg:augmentation} for a given $\mu$ and starting point
$\tilde x$. Then
$$
c^{\mu}(x,x_0) = c(x - x_0) - \mu \norm{x - x_0}_1 \leq 0,
$$
holds for all $x \in P \cap \Z^n$ and in particular for the optimal solution
$x^*$; this is simply the negation of the improvement condition. Now let $x'$ be any
iterate in the following call to Algorithm~\ref{alg:augmentation} for
which an augmentation step is performed with objective
$c^{\mu/2}(\cdot, \cdot)$ and starting point $x_0$, i.e., there exists
$x^+$ so that
$$
c^{\mu/2}(x^+,x') = c(x^+ - x') - \mu/2 \norm{x^+ - x'}_1 > 0.
$$
We can now combine these two inequalities, substituting \(x \leftarrow
x^*\), to obtain

\begin{align}
2 \frac{c(x^+ - x')}{\norm{x^+ - x'}_1} & > \mu \geq \frac{c(x^* - x_0)}{\norm{x^* - x_0}_1},
\end{align}
which implies
\begin{align}
\label{eq:geoScaleRestart}
c(x^+ - x') & \geq \frac{1}{2} \frac{\norm{x^+ - x'}_1}{\norm{x^* - x_0}_1} c(x^* - x_0) \geq \frac{1}{2n} c(x^* - x_0),
\end{align}
where $\norm{x^+ - x'}_1 \geq 1$ as the iterates are not identical and
$\norm{x^* - x_0}_1 \leq n$ as $x^*, x_0 \in P \subseteq [0,1]^n$. As
such each augmentation step recovers at least a
$\frac{1}{2n}$-fraction of the primal gap $c(x^* - x_0)$ and therefore
we can do at most $2n$ such iterations before the condition in
Line~\ref{eq:augStep} has to be violated. With this we can formulate
the geometric scaling restart scheme in
Algorithm~\ref{alg:geoScaling}. The analysis now is basically
identical to the one as for Algorithm~\ref{alg:bitScaling}, however
this time we have $O(\log n \norm{c}_\infty)$ restarts, leading to an
overall number of augmentation steps of $O(n \log n \norm{c}_\infty)$,
which can be further improved to $O(n \log \norm{c}_\infty)$, matching
that of bit-scaling, with the
simple observation in \cite{BPPP2015}.

\begin{algorithm}[htbp]
 \caption{\label{alg:geoScaling}Geometric Scaling} 
\begin{algorithmic}[1]
\REQUIRE Feasible solution \(x^0\)\\
\ENSURE Optimal solution of \(\max \set{cx}{x \in P \cap \Z^n}\)
    \STATE \(C \leftarrow \norm{c}_\infty + 1\), \(\mu \leftarrow nC\), \(\tilde{x}
    \leftarrow x^0\), \(c^{\mu}(x,y) \doteq c(x - y) - \mu \norm{x - y}_1. 
\) \COMMENT{initialization}
    \REPEAT
    \STATE \textbf{Call} \(\tilde x \leftarrow \mathcal A(\tilde x, c^\mu)\) \label{aeq:callStdAugGeo}
    \STATE \(\mu \leftarrow \mu/2\)
    \UNTIL{\(\mu < 1\)} \label{aeq:stopGeo}
  \RETURN \(\tilde{x}\) \COMMENT{return optimal solution}
  \end{algorithmic}
\end{algorithm}

\paragraph{Related Approaches.} Chvátal-Gomory cutting planes,
introduced by Chvátal in \cite{Chvatal}, are an important tool in
integer programming to approximate the integral hull
\(\conv(P \cap \Z^n)\) by means of successively strengthening an
initial relaxation \(P\) with \(\conv(P \cap \Z^n) \subseteq P\). This
is done by adding new inequalities valid for \(\conv(P \cap \Z^n)\)
cutting off chunks of \(P\) in each round.  A key question is how
many rounds of such strengthenings are needed until we recover
\(\conv(P \cap \Z^n)\). In \cite{CCH} it was shown that in general the
number of rounds can be arbitrarily large. It was then shown in
\cite{BEHS} via a restart argument that for the important case of
polytopes contained in \([0,1]^n\) the number of rounds can be upper
bounded by \(O(n^3 \log n)\). The key here is to use basic bounds on
the number of rounds, e.g., from \cite{CCH}, first for inequalities
with some maximum absolute entry \(c\), then doubling up \(c\) to
\(2c\), and restarting the argument. This bound was further improved
in \cite{ES} to \(O(n^2 \log n)\) by interleaving two
restart arguments, one multiplicative (e.g., doubling) over the
maximum absolute entry \(c\) and one additive (e.g., adding a
constant) over the dimension, which matches the lower bound of
\(\Omega(n^2)\) of \cite{rothvoss20170} up to a log factor; closing
this gap remains an open problem. As
mentioned in the context of the scheduled restarts of \cite{Roul17},
it might be possible that the additional log factor is due to the
restart schemes itself and removing it might require a different proof
altogether.

Another important application is the approximate Carathéodory problem,
where we want to approximate \(x^0 \in P\), where \(P\) is a polytope,
by means of a sparse convex combination \(x\) of vertices of \(P\), so
that \(\norm{x_0 - x} \leq \varepsilon\) for some norm
\(\norm{\cdot}\) and target accuracy \(\varepsilon\). In general it is
known that this can be done with a convex combination of
\(O(1/\varepsilon^2)\) vertices. However, it turns out as shown in
\cite{mirrokni17cara} that whenever \(x_0\) lies deep inside the
polytope \(P\), i.e., we can fit a ball around \(x_0\) with some
radius \(r\) into \(P\) as well, then we can exponentially improve
this bound via restarts to
\(O(\frac{1}{r^2} \log \frac{1}{\varepsilon})\). This restart argument
here is particularly nice. We run the original
\(O(1/\varepsilon^2)\)-algorithm down to some fixed accuracy and
obtain some approximation \(\tilde x\), then scale-up the feasible
region by a factor of \(2\), and restart the
\(O(1/\varepsilon^2)\)-algorithm on the residual \(x_0 - \tilde x\)
and repeat. The argument in \cite{mirrokni17cara} relies on mirror
descent as underlying optimization routine. More recently, it was
shown in \cite{cp2019approxCara} that the restarts can be removed and
adaptive bounds for more complex cases can be obtained by
using conditional gradients as base optimization algorithm, which
automatically adapts to sharpness (and optima in the interior)
\cite{YiAdapativeFW,KDP2018}.

\section{Submodular Function Maximization}
\label{sec:sfm}

We now turn our attention to submodular function
maximization. Submodularity captures the \emph{diminishing returns}
property and is widely used in optimization and machine learning. In
particular, we will consider the basic but important setup of maximizing a monotone,
non-negative, submodular function subject to a single cardinality
constraint. To this end we will briefly repeat necessary notions.  A
set function $g:2^V\to\R_+$ is \emph{submodular} if and only if for
any $e\in V$ and $A\subseteq B\subseteq V\backslash\{e\}$ we have
$g_A(e)\geq g_B(e)$, where $g_A(e) \doteq g(A+e)-g(A)$ denotes the
\emph{marginal gain of \(e\) w.r.t. \(A\)} and
$A+e \doteq A\cup\{e\}$, slightly abusing notation. The submodular
function $g$ is \emph{monotone} if for all $A\subseteq B\subseteq V$
it holds $g(A)\leq g(B)$ and \emph{non-negative} if \(g(A) \geq 0\)
for all \(A \subseteq V\).

Given a monotone, non-negative submodular
function \(g\) over ground set \(V\) of size \(n\) and a budget \(k\), we consider
the problem

\begin{equation}
  \label{eq:sfm}
  \max_{S \subseteq V, \card{S} \leq k} g(S)
\end{equation}

It is well known that solving Problem~\eqref{eq:sfm} exactly is
NP-hard under the value oracle model, however the greedy algorithm
(Algorithm~\ref{alg:greedy}) that in each iteration adds the element
that maximizes the marginal gain yields a \((1-1/\mathrm{e})\)-approximate
solution \(S^+ \subseteq V\) with \(\card{S^+} \leq k\),
i.e., \(g(S^+) \geq (1-1/\mathrm{e})\ g(S^*)\), where \(S^*
=\argmax_{S\subseteq V, \card{S}\leq k} g(S) \subseteq V\) is an
optimal solution to Problem~\eqref{eq:sfm} and \(\mathrm{e}\) denotes
the Euler constant (see
\cite{nemhauser1978analysis,fisher1978analysis}).

\begin{algorithm}[htbp]
\caption{Greedy Algorithm}\label{alg:greedy}
\begin{algorithmic}[1]
  \REQUIRE Ground set $V$ of size \(n\), budget \(k\), and monotone,
  non-negative, submodular function $g:2^{V} \rightarrow \R_+$.
  \ENSURE feasible set $S^+$ with \(\card{S^+} \leq k\).  \STATE
  $S^+\leftarrow\emptyset$ \WHILE{\(\card{S^+} \leq k\)}
  \STATE
  \label{aeq:oracle}\(e \leftarrow \argmax_{e \in V \setminus S^+}
  g_{S^+}(e)\) \STATE \(S^+ \leftarrow S^+ + e\) \ENDWHILE
\end{algorithmic}
\end{algorithm}

The proof of the approximation guarantee of \(1-1/\mathrm{e}\) is based on the
insight that in each iteration it holds:
\begin{equation}
  \label{eq:factorSFM}
g(S^{*})-g(S^+) \leq k \cdot \max_{e\in V}g_{S^+}(e).
\end{equation}
To see that Inequality~\eqref{eq:factorSFM} holds, let
$S^{*}=\{e_{1},\dots,e_{k}\}$, then
\begin{align*}
g(S^{*}) & \leq g(S^{*} \cup S^+) =g(S^+)+\sum_{i=1}^k g_{S^+\cup\{e_{1},\dots,e_{i-1}\}}(e_i)\\
 & \leq g(S^+)+\sum_{i=1}^k g_{S^+}(e_i) \leq g(S^+) + k\max_{e\in V}g_{S^+}(e),
\end{align*}
where the first inequality follows from monotonicity, the equation
follows from the definition of $g_S(v)$, the second inequality from submodularity, and the last inequality from
taking the maximizer. 

With Inequality~\eqref{eq:factorSFM} the proof of the \((1-1/\mathrm{e})\)-approximation is immediate. In each
iteration the greedy element we add satisfies \(\max_{e\in
  V}g_{S^+}(e) \geq \frac{1}{k} (g(S^{*})-g(S^+))\), therefore after
\(k\) iterations we have obtained a set \(S^+\) with \(\card{S^+} =
k\), with
\[g(S^*) - g(S^+) \leq (1-1/k)^k (g(S^*) - g(\emptyset)) \leq
  (1-1/k)^k g(S^*) \leq \frac{1}{\mathrm{e}} g(S^*),\]
so that the desired guarantee \((1-\frac{1}{\mathrm{e}}) g(S^*) \leq   g(S^+)\) follows.

Unfortunately, due to Line~\ref{aeq:oracle} in
Algorithm~\ref{alg:greedy} computing such a \((1-1/\mathrm{e})\)-approximate
solution can cost up to \(O(k n)\) evaluations of \(g\) in the value
oracle model, where we can only query function values of \(g\). For
realistic functions this is often quite prohibitive. We will now see a
different application of a restart scheme to reduce the total number of
function evaluations of \(g\) by allowing for a small error
\(\varepsilon > 0\). We obtain a total number of evaluations of
\(g\) of \(O(\frac{n}{\varepsilon} \log \frac{n}{\varepsilon})\),
quasi-linear and independent of \(k\), to compute a
\((1-1/\mathrm{e}-\varepsilon)\)-approximate solution. The argument is due to
\cite{badanidiyuru2014fast} and similar in nature to the argument in
Section~\ref{sec:discopt} for geometric scaling. We simplify the argument slightly for
exposition; see \cite{badanidiyuru2014fast} for details. 

The basic idea is rather than computing the actual maximum in Line~\ref{aeq:oracle} in
Algorithm~\ref{alg:greedy}, we collect all elements of marginal gains that are
roughly maximal within a \((1-\varepsilon)\)-factor, then scale down
the estimation of the maximum, and then restart. We present the
restart scheme, the so-called Threshold Greedy Algorithm in
Algorithm~\ref{alg:TGreedy}. This time we present the scheme and the
base algorithm directly together. Note that the inner loop in
Lines~\ref{aeq:innerStart} to \ref{aeq:innerEnd} in
Algorithm~\ref{alg:TGreedy} adds all elements that have approximately
maximal marginal gain. The restarts are happening whenever we go back to the
beginning of the outer loop starting in Line~\ref{aeq:restartsTG},
with a reset value for \(\Phi\).

\begin{algorithm}[htbp]
\caption{Threshold Greedy Algorithm}\label{alg:TGreedy}
\begin{algorithmic}[1]
\REQUIRE Ground set $V$ of size \(n\), budget \(k\), accuracy \(\varepsilon\), and monotone,
non-negative, 
submodular function $g:2^{V} \rightarrow \R_+$
\ENSURE feasible set $S^+$ with \(\card{S^+} \leq k\).
\STATE $S^+\leftarrow\emptyset$, \(\Phi_0 \leftarrow \max_{e \in V}
g(e)\), \(\Phi \leftarrow \Phi_0\)
\WHILE{\(\Phi \geq \frac{\varepsilon}{n}\Phi_0\)} \label{aeq:restartsTG}
\FOR{\(e \in V\)}\label{aeq:innerStart}
\IF{\(\card{S^+} < k\) and \(g_{S^+}(e) \geq \Phi\)} \label{aeq:appxMax}
\STATE \(S^+ \leftarrow S^+ + e\) \label{aeq:addEle}
\ENDIF
\ENDFOR \label{aeq:innerEnd}
\STATE \(\Phi
  \leftarrow \Phi (1- \varepsilon)\)
\ENDWHILE 

\end{algorithmic}
\end{algorithm}

We will first show that the gain from any new element \(e \in V\) added in
Line~\ref{aeq:addEle} of Algorithm~\ref{alg:TGreedy} is at least

\begin{equation}
  \label{eq:gain}
    g_{s^+}(e) \geq \frac{1-\varepsilon}{k} \sum_{x \in S^* \setminus S^+ } g_{S^+}(x).
\end{equation}
To this end suppose we have have chosen element \(e \in V\) to be
added. Then \(g_{s^+}(e) \geq \Phi\) by Line~\ref{aeq:appxMax} and for
all \(x \in S^* \setminus (S^+ + e)\) we have \(g_{S^+}(x) \leq
\Phi/(1-\varepsilon)\); otherwise we would have added \(x\) in an earlier
restart with a higher value \(\Phi\) already. Combining the two
inequalities we obtain
\begin{equation}
  \label{eq:approxGain}
  g_{s^+}(e) \geq (1-\varepsilon) g_{S^+}(x),
\end{equation}
for all  \(x \in S^* \setminus (S^+ + e)\) and averaging those
inequalities leads to
\begin{equation}
  \label{eq:approxGainAvg}
  g_{s^+}(e) \geq \frac{1-\varepsilon}{\card{S^* \setminus S^+}}
  \sum_{x \in S^* \setminus S^+} g_{S^+}(x) \geq \frac{1-\varepsilon}{k}
  \sum_{x \in S^* \setminus S^+}  g_{S^+}(x),
\end{equation}
which is the desired inequality. From this we immediately recover the
(approximate) analog of Inequality~\eqref{eq:factorSFM}. We have via
submodularity and non-negativity
\begin{equation}
  \sum_{x \in S^* \setminus S^+}  g_{S^+}(x) \geq g_{S^+}(S^*) \geq
  g(S^*) - g(S^+),
\end{equation}
and together with Inequality~\eqref{eq:approxGainAvg}
\begin{equation}
  \label{eq:factorSFMApprox}
  g_{s^+}(e) \geq \frac{1-\varepsilon}{k} (g(S^*) - g(S^+)).
\end{equation}

Therefore, as before, after
\(k\) iterations we obtain a set \(S^+\) with \(\card{S^+} =
k\), with
\begin{align*}
  g(S^*) - g(S^+) & \leq (1-(1-\varepsilon)/k)^k (g(S^*) - g(\emptyset)) \leq
                    (1-(1-\varepsilon)/k)^k g(S^*) \\
  & \leq \frac{1}{\mathrm{e}^{(1-\varepsilon)}} g(S^*) \leq
    \left( \frac{1}{\mathrm{e}} + \varepsilon \right ) g(S^*),
\end{align*}
leading to our guarantee \(g(S^+) \geq \left( 1 - \frac{1}{\mathrm{e}}
    -\varepsilon \right) g(S^*)\). If we do fewer than \(k\)
  iterations, the total gain of all remaining elements is less than
  \(\varepsilon\), establishing the guarantee in that case. 

 Now for the number of evaluations of \(g\), first consider the loop in
 Line~\ref{aeq:restartsTG} of Algorithm~\ref{alg:TGreedy}. The loops
 stops after \(\ell\) iterations, whenever \((1-\varepsilon)^\ell \leq
 \frac{\varepsilon}{n}\), which is satisfied if
 \(1/(1-\varepsilon)^\ell \geq (1+\varepsilon)^\ell \geq
 \frac {n}{\varepsilon}\) and hence \(\ell \geq \frac{1}{\varepsilon}
 \log \frac{n}{\varepsilon}\). For each such loop iteration we have
 at most \(O(n)\) evaluations of \(g\) in Line~\ref{aeq:appxMax}, leading
 to the overall bound of \(O(\frac{n}{\varepsilon}
 \log \frac{n}{\varepsilon})\) evaluations of \(g\). 

\paragraph{Related Approaches.}  

The approach presented here for the basic case with a single
cardinality constraint can be applied more widely as already done in
\cite{badanidiyuru2014fast} for matroid, knapsack, and \(p\)-system
constraints. It can be also used to reduce the number of evaluations
in the context of robust submodular function maximization
\cite{submodular2017,submodular2017jour}.

A similar restart approach has been used to \lq{}lazify\rq{}
conditional gradient algorithms in
\cite{LPZZ2017,BPZ2017,BPZ2017jour}. Here is the number of calls to
the underlying linear optimization oracle is dramatically reduced by
reusing information from previous iterations by solving the linear
optimization problem only approximately as done in the case of the
Threshold Greedy Algorithm. The algorithm, in
a similar vein, is then restarted, whenever the threshold for
approximation of the maximum is too large.

\section*{Acknowledgement}
\label{sec:acknowledgement}

We would like to thank Gábor Braun and Marc Pfetsch for helpful
comments and feedback on an earlier version of this article. 

\bibliographystyle{splncs04}
\bibliography{references}

\end{document}